\documentclass[english,envcountsect, envcountsame]{svjour3}
\usepackage[T1]{fontenc}
\usepackage[latin9]{inputenc}
\usepackage{bbding}
\usepackage{mathrsfs}
\usepackage{url}
\usepackage{amsmath}
\usepackage{amssymb}
\PassOptionsToPackage{normalem}{ulem}
\usepackage{ulem}

\makeatletter
\newenvironment{svmultproof}{\begin{proof}}{\qed\end{proof}}

\usepackage{amsmath,amsfonts}
\usepackage{enumitem}
\DeclareMathOperator{\sa}{sa}

\DeclareMathOperator{\At}{At}
\setlist[enumerate,1]{label = (\alph*)}

\makeatother

\usepackage{babel}

\begin{document}
\title{Countably-categorical Boolean rings with distinguished ideals}
\author{Andrew B. Apps\thanks{I gratefully acknowledge the use of the library facilities of the
University of Cambridge.\protect\phantom{\Envelope{}}}}
\institute{\Envelope{}  Andrew B. Apps, Independent researcher, St Albans, UK\\
$^{\phantom{\Envelope{}}}$ andrew.apps@apps27.co.uk\\\\}
\date{August 2025}
\maketitle
\begin{abstract}
We describe and classify countable Boolean rings (which may or may
not have a multiplicative identity) with finitely many distinguished
ideals whose elementary theory is countably categorical. This extends
the description by Macintyre and Rosenstein and subsequent authors
of countably categorical Boolean algebras with finitely many distinguished
ideals. Following Pierce, we take a topological approach using the
language of PO systems (partially ordered sets with a distinguished
subset) and topological Boolean algebras. We provide two different
classifications via invariants that uniquely determine the isomorphism
type: one using finite PO systems and the other using finite posets.We
discuss how our findings link with previous results, but the paper
is otherwise self-contained.
\end{abstract}

\section{Introduction}

A countable structure $M$ for a language $L$ is \emph{$\omega$-categorical
}if $M$ is determined up to isomorphism within the class of countable
$L$-structures by its first order properties. By the theorem of Engeler,
Ryll-Nardzewski and Svenonius~\cite{Engeler,Ryll,Svenonius}, this
is equivalent to saying that for all $r$, the automorphism group
of $M$ is \emph{almost $r$-transitive }on $M$: that is, it has
only finitely many orbits on $M^{r}$.

The problem of describing $\omega$-categorical Boolean algebras with
finitely many distinguished ideals arose during the classification
by Macintyre and Rosenstein~\cite{MacRos} of $\omega$-categorical
rings without nilpotent elements. If $C_{1},\ldots,C_{n}$ are closed
subsets of a compact Stone space $X$ and $J_{i}$ is the ideal of
$R$ (the Boolean algebra underlying $X$) corresponding to $C_{i}$
under the Stone correspondence, they showed that the Boolean algebra
system $(R,J_{1},\ldots,J_{n})$ with distinguished ideals is $\omega$-categorical
iff $\langle J_{1},\ldots,J_{n}\rangle$ is a finite sub-algebra of
the set $H(R)$ of all ideals of $R$, viewed as a Heyting algebra,
and $R/J$ has only finitely many atoms for each $J\in\langle J_{1},\ldots,J_{n}\rangle$.

There have been a number of subsequent studies of this condition (see~\cite{Alaev,AppsBP,Palchunov,Touraille}),
mostly from the viewpoint of Heyting algebras. In this paper we adopt
a topological approach that looks at partitions of the Stone space,
working with partially ordered sets rather than the more complex structure
of Heyting algebras. This is a natural approach, as Macintyre and
Rosenstein's result originally arose as a question about distinguished
closed subsets of a Stone space. We generalise by considering Boolean
rings which may not have an identity rather than just Boolean algebras,
as both compact and locally compact Stone spaces arise naturally in
the study of $\omega$-categorical structures. 

Finally, we classify $\omega$-categorical Boolean rings with finitely
many distinguished ideals by identifying relevant invariants, firstly
in terms of finite PO systems (posets with a distinguished subset)
and secondly in terms of finite posets. By~\cite[Theorem~5]{MacRos},
Boolean structures of this form can be combined with a finite structure
$M$ and some distinguished substructures of $M$ to create \emph{``filtered
Boolean extensions''} which are also $\omega$-categorical; Schmerl~\cite{Schmerl}
extended this result to include infinite structures $M$ where $M$
together with the distinguished substructures is $\omega$-categorical.
It is therefore of some interest to understand and classify the Boolean
structures that underpin these filtered Boolean extensions.

First, we recall some definitions.
\begin{definition}
If $X$ is a topological space and $S\subseteq X$, we write $S^{\prime}$
for the \emph{derived set }of $S$, namely the set of all limit points
of $S$, and $\overline{S}$ for the closure of $S$. We say that
$S$ is \emph{crowded} (or \emph{dense-in-itself}) if it has no isolated
points; equivalently, if $S\subseteq S^{\prime}$.

A \emph{topological Boolean Algebra (TBA)} is a Boolean algebra with
an additional unary operator $^{\prime}$ satisfying $0^{\prime}=0$,
$(A\cup B){}^{\prime}=A^{\prime}\cup B^{\prime}$ and $A^{\prime}{}^{\prime}\subseteq A{}^{\prime}$;
and a \emph{closure algebra }is a Boolean algebra together with an
additional closure operator $^{-}$ such that $A\subseteq\overline{A}$,
$\overline{\overline{A}}=\overline{A}$, $\overline{A\cup B}=\overline{A}\cup\overline{B}$,
and $\overline{0}=0$. A TBA ($\mathscr{D},^{\prime})$ can also be
viewed as a closure algebra ($\mathscr{D},^{-})$ by setting $\overline{B}=B\cup B^{\prime}$
for $B\in\mathscr{D}$. If $X$ is a topological space, then $(2^{X},^{-})$
and $(2^{X},^{\prime})$ become a closure algebra and a TBA respectively.
For $A\subseteq2^{X}$, we write $\langle A,^{-}\rangle$ and $\langle A,^{\prime}\rangle$
for the closure sub-algebra and sub-TBA of $2^{X}$ respectively generated
by $A$, which each include $X$ itself.

We will say that a partition $\mathscr{X}$ of a topological space
$X$ is \emph{finite-crowded }if $\mathscr{X}$ is finite and every
element of $\mathscr{X}$ is either finite or crowded.
\end{definition}
\emph{Notation}: if $(J_{1},\ldots,J_{n})$ is an $n$-tuple, we will
write $\boldsymbol{J}^{(n)}$ for $(J_{1},\ldots,J_{n})$ (i.e.\ when
the order matters) and $\boldsymbol{J}^{\{n\}}$ for $\{J_{1},\ldots,J_{n}\}$,
and similarly for $\boldsymbol{C}^{(n)}$, $\boldsymbol{Q}^{\{n\}}$
etc. 

If $C$ is a closed subset of the Stone space $X$ of the Boolean
ring $R$, we write $I(C)$ for the corresponding ideal of $R$, namely
$I(C)=\{A\in R\mid A\cap C=\emptyset\}$. For a subset $\mathscr{E}$
of $2^{X}$, we let $(R,I(\overline{\mathscr{E}}))$ denote the Boolean
ring $R$ with distinguished ideals $\{I(\overline{A})\mid A\in\mathscr{E}\}$. 

Our primary result is the following, which extends the well-known
result of Macintyre and Rosenstein~\cite[Theorem~7]{MacRos}. 
\begin{description}
\item [{\textbf{Theorem~\ref{Theorem: main result}}}] \emph{(partial
statement) Let $\boldsymbol{J}^{(n)}$ be distinguished ideals of
the countable Boolean ring $R$, corresponding to closed subsets $\boldsymbol{C}^{(n)}$
of the Stone space $X$ of $R$. Then the following are equivalent:}
\begin{enumerate}
\item \emph{$(R,\boldsymbol{J}^{(n)})$ is $\omega$-categorical; }
\item \emph{$\langle\boldsymbol{C}^{\{n\}},^{-}\rangle$ is finite, and
each of its atoms has only finitely many isolated points;}
\item \emph{$\langle\boldsymbol{C}^{\{n\}},^{\prime}\rangle$ is finite,
and its atoms form a finite-crowded partition of $X$.}
\end{enumerate}
\end{description}
We will make frequent use of the duality established by Pierce~\cite{Pierce}
between finite TBAs generated by their closed elements and finite
\emph{PO systems:}
\begin{definition}
A \emph{PO system }is a set $P$ with an anti-symmetric transitive
relation $<$; equivalently, it is a poset with a distinguished subset
$P_{1}$, where $p<p$ iff $p\in P_{1}$. 

If $\mathscr{D}$ is a TBA or closure algebra, we let $\At(\mathscr{D})$
denote the set of atoms of $\mathscr{D}$ and $1_{\mathscr{D}}$ denote
the largest element of $\mathscr{D}$. 

If $\mathscr{D}$ is a TBA and $P$ is a finite PO system, we will
say that a partition $\mathscr{X}$ of $1_{\mathscr{D}}$ is a \emph{complete
$P$-partition of $1_{\mathscr{D}}$ }if there is a labelling $\mathscr{X}=\{X_{p}\mid p\in P\}$
such that $\ensuremath{X_{p}^{\prime}=\bigcup_{q<p}X_{q}\text{ for all }p\in P}$.
We will relax the requirement for $P$ to be finite if $\mathscr{D}$
is a sub-TBA of $2^{X}$ for a topological space $X$, as infinite
unions exist in this case.
\end{definition}
Under this duality, if $X$ is a Stone space and $\mathscr{D}$ is
a sub-TBA of $2^{X}$ such that $\At(\mathscr{D})$ forms a finite-crowded
partition of $X$, then $\At(\mathscr{D})$ becomes a complete $P$-partition
of $X$ for a suitable PO system $P$. Our next result establishes
$\omega$-categoricity in this situation and is the key step in showing
that~\ref{enu:Thm4.3} implies~\ref{enu:Thm4.1} in Theorem~\ref{Theorem: main result}:
\begin{description}
\item [{\textbf{Theorem~\ref{Thm: P-partition categorical}}}] \emph{Let
$P$ be a PO system and $\mathscr{X}$ a complete $P$-partition of
the Stone space of the countable Boolean ring $R$. Then the system
$(R,I(\overline{\mathscr{X}}))$ is $\omega$-categorical iff $\mathscr{X}$
is finite-crowded.}
\end{description}
These Theorems will lead in turn to a classification of $\omega$-categorical
Boolean rings with distinguished ideals firstly in terms of finite
PO systems (Section~\ref{sec:Classification-POsystem}) and secondly
in terms of finite posets (Section~\ref{sec:Classificationposet}). 

With the exception of Section~\ref{sec:Links-with-previous}, the
paper is self-contained. For example, the topological methods used
here provide a different proof of the original result of Macintyre
and Rosenstein~\cite[Theorem~7]{MacRos}. They also extend that result
and others from Boolean rings with a $1$ to Boolean rings without
a $1$, corresponding to compact and locally compact Stone spaces
respectively. Where this is the case, we provide a brief commentary
after the relevant Theorem, and a more detailed explanation of the
inter-relationships in Section~\ref{sec:Links-with-previous}, including
the translation between categories where existing results are stated
in terms of Heyting algebras. There are two parallel sets of category
dualities at play here: firstly that between finite closure algebras
generated by their closed elements, finite posets, and finite Heyting
algebras; and secondly that between finite TBAs generated by their
closed elements, finite PO systems, and finite Heyting ``sa-algebras''. 

\subsection{Notation}

If $A$ and $B$ are sets, we write $A\subseteq B$ to denote that
$A$ is contained in $B$, $A-B$ for $\{x\in A\mid x\notin B\}$,
$A\dotplus B$ for the disjoint union of $A$ and $B$, and $|A|$
for the cardinal number of $A$. We write $\mathbb{N}$ and $\mathbb{N}_{+}$
for the non-negative and positive integers respectively. 

If $R$ is a Boolean ring with Stone space $X$, we identify $R$
with the ring of compact open subsets of $X$ under the Stone correspondence,
whereby atoms of the Boolean ring correspond to isolated points of
the Stone space. $R$ has a $1$ iff $X$ is compact, and is countable
iff $X$ is second countable (i.e.\ the topology has a countable
base of open sets); we use the term \emph{$\omega$-Stone space }to
denote the Stone space of a countable Boolean ring. For $A\in R$,
we write $(A)$ for the ideal $\{B\in R\mid B\subseteq A\}$. 

\section{\label{sec:Prior-results}Complete finite $P$-partitions of a Stone
space}

Let $P$ be a finite PO system and let $\mathscr{X}=\{X_{p}\mid p\in P\}$
be a complete $P$-partition of the Stone space $X$ of the Boolean
ring $R$. For $A\in R$, we will define a ``measure'' $\mu(A)$,
being a formal sum of the form $\sum_{p\in F}n_{p}.p$, where $F\subseteq P$
and $n_{p}\in\mathbb{N}_{+}$, and we shall see that the pair $(A,\mathscr{X}_{A})$,
where $\mathscr{X}_{A}$ is the restriction of $\mathscr{X}$ to $A$,
is determined up to homeomorphism by $\mu(A)$. This will be the crucial
step in establishing the ``if'' statement in Theorem~\ref{Thm: P-partition categorical}.

The following definitions and notation are adapted from~\cite{Apps-Stone}.
\begin{definition}
For $p,q\in P$, we write $p\leqslant q$ to mean $p<q$ or $p=q$;
$P^{d}$ for $\{p\in P\mid p\nless p\}$; and $P_{\min}$ for the
minimal elements of $(P,\leqslant)$. A subset $L$ of $P$ is a \emph{lower
}(respectively\emph{ upper})\emph{ subset }of $P$ if $p\in L$ and
$q<p$ (respectively $q>p$) implies $q\in L$.

For a subset $A$ of $X$, we define its \emph{type }$T(A)=\{p\in P\mid A\cap X_{p}\neq\emptyset\}$. 

An element $A\in R$ is \emph{$p$-trim,} and we write $t(A)=p$,
if $T(A)=\{q\in P\mid q\geqslant p\}$, with also $|A\cap X_{p}|=1$
if $p\in P^{d}$; and is \emph{trim }if it is $p$-trim for some $p\in P$. 

For $A\in R$, the \emph{$\mathscr{X}$-measure of $A$} is the formal
sum $\mu(A)=\sum_{p\in F}n_{p}.p$, where $F=T(A)_{\min}$, $n_{p}=1$
for $p\in F-P^{d}$ and $n_{p}=|A\cap X_{p}|$ for $p\in F\cap P^{d}$.
A partition of $A$ into trim sets is a \emph{$\mu(A)$-partition
}if it contains precisely $n_{p}$ $p$-trim sets for each $p\in F$;
as we shall see, this equates to a minimum decomposition of $A$ into
disjoint trim sets. We set $\mu(0)=0$.

If $\{A_{i}\mid i\in I\}$ are subsets of $X$, we say that a partition
$\{X_{p}\mid p\in P\}$ of $X$ \emph{refines }$\{A_{i}\mid i\in I\}$
if for each $i\in I$, we can find $Q_{i}\subseteq P$ such that $A_{i}=\bigcup_{q\in Q_{i}}X_{q}$.

If $\{X_{p}\mid p\in P\}$ and $\{Y_{p}\mid p\in P\}$ are complete
$P$-partitions of the Stone spaces $X$ and $Y$ respectively, a
homeomorphism $\alpha\colon X\rightarrow Y$ is a \emph{$P$-homeomorphism}
if $X_{p}\alpha=Y_{p}$ for all $p\in P$.
\end{definition}
\begin{remark}
\label{rem:trimpartitions}The broader definitions of ``trim'' and
``semitrim'' partitions in~\cite{Apps-Stone} apply where $P$
is an infinite PO system or poset and where $\{X_{p}\mid p\in P\}$
is a partition of a dense subset of $X$ rather than a complete partition
of the whole of $X$. 

If $\{X_{p}\mid p\in P\}$ is a complete $P$-partition of $X$ for
a PO system $P$, then it follows that $\ensuremath{\overline{X_{p}}=\bigcup_{q\leqslant p}X_{q}\text{ for all }p\in P}$,
as $\overline{X_{p}}=X_{p}\cup X_{p}^{\prime}$; that $X_{p}$ has
no isolated points if $p\notin P^{d}$; and that $X_{p}$ is discrete
if $p\in P^{d}$.

For $A\in R$ and $p\in P^{d}\cap T(A)_{\min}$, consideration of
accumulation points shows that $|A\cap X_{p}|$ is necessarily finite.
Hence $X_{p}$ is finite iff $X_{p}$ is compact and $p\in P_{\min}^{d}$. 
\end{remark}
We will need the following result of Apps~\cite{Apps-Stone}, noting
that for finite $P$, a complete $P$-partition is always a ``trim''
partition (as defined in~\cite{Apps-Stone}). We include a fresh
proof which is simpler than that in~\cite{Apps-Stone}, in order
to make this paper as self-contained as possible. 
\begin{theorem}
\label{Thm:same mu=00003Dhomeom}(\cite[Theorem~3.5]{Apps-Stone})
Let $P$ be a finite PO system and $\mathscr{X}$ and $\mathscr{Y}$
complete $P$-partitions of the Stone spaces $X$ and $Y$ of countable
Boolean rings $R$ and $S$ respectively. Suppose $A\in R$ and $B\in S$
are such that $\mu(A)=\mu(B)$. Then there is a $P$-homeomorphism
$\alpha\colon A\rightarrow B$.
\end{theorem}

\subsection{Proof of Theorem~\ref{Thm:same mu=00003Dhomeom}}

We first recall some basic properties of complete $P$-partitions,
where $P$ is a finite PO system.
\begin{proposition}
\label{Propn std split of trim set}Let $P$ be a finite PO system,
let $\{X_{p}\mid p\in P\}$ be a complete $P$-partition of the Stone
space $X$ of the Boolean ring $R$ and let $A\in R$. Then:
\begin{enumerate}
\item \label{enu:Prop1.1}\cite[Proposition~2.12(iii)]{Apps-Stone} $x\in X_{p}$
iff $x$ has a neighbourhood basis of $p$-trim sets;
\item \label{enu:Prop1.2} if $A$ is $p$-trim and $q>p$, then we can
write $A=B\dotplus C$, where $B$ is $p$-trim and $C$ is $q$-trim;
\item \label{enu:Prop1.3}there is a partition $\mathscr{B}$ of $A$ into
disjoint sets such that each $B\in\mathscr{B}$ is a compact open
subset of a trim set; and any such partition $\mathscr{B}$ has a
subset $\mathscr{C}$ such that $\mathscr{C}$ is a $\mu(A)$-partition
of $\bigcup\mathscr{C}$;
\item \label{enu:Prop1.4}\cite[Proposition~3.2]{Apps-Stone} $A$ has a
$\mu(A)$-partition.
\end{enumerate}
\end{proposition}
\begin{svmultproof}
\ref{enu:Prop1.1} Suppose $B\in R$ contains $x\in X_{p}$; we need
to find a $p$-trim neighbourhood of $x$ contained in~$B$. If $p\in P^{d}$,
reduce $B$ if necessary to assume that $B\cap X_{p}=\{x\}$. For
each $q\ngeqslant p$ we have $x\notin\overline{X_{q}}$, so we can
find $C_{q}\in R$ with $x\in C_{q}\subseteq B$ and $C_{q}\cap\overline{X_{q}}=0$.
Then $\bigcap_{q\ngeqslant p}C_{q}$ is the required $p$-trim neighbourhood
of $x$, as if $q\geqslant p$ then all neighbourhoods of $x$ meet
$X_{q}$. The converse follows easily, as if $x\in X_{q}$ $(q\neq p)$,
then $x$ cannot also have a neighbourhood base of $p$-trim sets.

\ref{enu:Prop1.2} Choose $x\in A\cap X_{q}$ and a $q$-trim $C\subseteq A$
such that $x\in C$, ensuring that $(A-C)\cap X_{p}\neq\emptyset$
if $q=p$ (when $X_{p}$ is crowded). Let $B=A-C$, which is $p$-trim.

\ref{enu:Prop1.3} For each $y\in A$, use~\ref{enu:Prop1.1} to
find a trim set $D_{y}\subseteq A$ containing~$y$. By compactness,
we can find a finite subset $\{D_{y_{1}},\ldots,D_{y_{n}}\}$ of trim
sets that cover~$A$. Let $B_{j}=D_{y_{j}}-\bigcup_{k<j}D_{y_{k}}$
for each $j\leqslant n$; then $\mathscr{B}=\{B_{j}\mid j\leqslant n\}-\{\emptyset\}$
has the required properties. Let $u(B_{j})=t(D_{y_{j}})$.

Now let $F=T(A)_{\min}$, and for $p\in F$ let $\mathscr{B}_{p}=\{B\in\mathscr{B}\mid p\in T(B)\}$,
which is non-empty. If $B\in\mathscr{B}_{p}$, then $u(B)=p$ as $p\in T(A)_{\min}$,
and so $B$ is $p$-trim; hence also $\mathscr{B}_{p}\cap\mathscr{B}_{q}=\emptyset$
for $p\neq q$. Let $\mathscr{C}_{1}$ contain one element of $\mathscr{B}_{p}$
for each $p\in F-P^{d}$, and let $\mathscr{C}=\mathscr{C}_{1}\cup\bigcup\{\mathscr{B}_{p}\mid p\in F\cap P^{d}\}$
to finish.

\ref{enu:Prop1.4} Let $\mathscr{B}$ and $\mathscr{C}$ be as in~\ref{enu:Prop1.3},
and let $E=\bigcup\mathscr{C}$, so that $\mathscr{C}$ is a $\mu(A)$-partition
of $E$. Let $\mathscr{D}=\mathscr{B}-\mathscr{C}$. For each $D\in\mathscr{D}$,
we can find $r_{D}\in F$ with $r_{D}\leqslant u(D)$ and $C_{D}\in\mathscr{C}$
such that $t(C_{D})=r_{D}$. Now $D\cup C_{D}$ is still $r_{D}$-trim,
as $T(D)\subseteq T(C_{D})$, and if $r_{D}\in F\cap P^{d}$ then
$D\cap X_{r_{D}}=\emptyset$ by construction of $\mathscr{C}$. So
we can replace $C_{D}$ with $C_{D}\cup D$ to obtain a $\mu(A)$-partition
of $E\cup D$. Proceeding in this way, we obtain a $\mu(A)$-partition
of~$A$.
\end{svmultproof}

To prove Theorem~\ref{Thm:same mu=00003Dhomeom}, we will use the
following version of Vaught's Theorem, as cited in~\cite[1.1.3]{PierceMonk}.
\begin{theorem}
\label{(Vaught)}(Vaught) Suppose $\sim$ is a relation between elements
of countable Boolean algebras $R$ and $S$ such that:
\begin{enumerate}
\item $1_{R}\sim1_{S}$;
\item if $A\sim0_{S}$, then $A=0_{R}$; and vice versa;
\item \label{enu: Vaught3}if $A\sim(B_{1}\dotplus B_{2})$ ($B_{1},B_{2}\in S$),
then we can write $A=A_{1}\dotplus A_{2}$, with $A_{i}\in R$ and
$A_{i}\sim B_{i}$ $(i=1,2)$; and vice versa.
\end{enumerate}
Then there is an isomorphism $\alpha\colon R\rightarrow S$ such that
each $A\in R$ can be expressed as $A=A_{1}\dotplus\ldots\dotplus A_{n}$
where $A_{i}\sim A_{i}\alpha$ for all $i\leqslant n$.
\end{theorem}
\begin{svmultproof}
\emph{of Theorem~\ref{Thm:same mu=00003Dhomeom}}

Let $P$ be a finite PO system and $\mathscr{X}=\{X_{p}\mid p\in P\}$
and $\mathscr{Y}=\{Y_{p}\mid p\in P\}$ be complete $P$-partitions
of the Stone spaces $X$ and $Y$ of countable Boolean rings $R$
and $S$ respectively. Suppose $A\in R$ and $B\in S$ are such that
$\mu(A)=\mu(B)$. We must find a $P$-homeomorphism between $A$ and
$B$. 

Define a relation between the Boolean algebras $(A)$ and $(B)$ by
setting $C\sim D$ iff $\mu(C)=\mu(D)$, for $C\subseteq A$ and $D\subseteq B$.
Clearly $C\sim0$ iff $C=0$, and $A\sim B$. 

We must show that $\sim$ satisfies criterion~\ref{enu: Vaught3}
of Theorem~\ref{(Vaught)}. By symmetry, it is enough to consider
the case where $C\subseteq A$, $C\sim D$ and $D=D_{1}\dotplus D_{2}$,
where $D_{1},D_{2}\subseteq B$. If $\mathscr{C}$ and $\mathscr{D}$
are partitions of subsets of $A$ and $B$ respectively with each
partition element trim, we write $\mathscr{C}\sim\mathscr{D}$ to
mean that there is a type-preserving bijection from $\mathscr{C}$
to $\mathscr{D}$. 

Let $\mathscr{B}_{i}$ be a $\mu(D_{i})$-partition of $D_{i}$ ($i=1,2)$
and let $\mathscr{A}$ be a $\mu(C)$-partition of~$C$. By Proposition~\ref{Propn std split of trim set}\ref{enu:Prop1.3},
$\mathscr{B}_{1}\cup\mathscr{B}_{2}$ must contain a subset $\mathscr{B}_{3}$,
say, such that $\mathscr{A}\sim\mathscr{B}_{3}$; let $\mathscr{B}_{4}=(\mathscr{B}_{1}\cup\mathscr{B}_{2})-\mathscr{B}_{3}$.
For each $E\in\mathscr{B}_{4}$, we can find $F\in\mathscr{A}$ such
that $t(F)\leqslant t(E)$; but (counting) we cannot have $t(F)=t(E)\in P^{d}$,
and so $t(F)<t(E)$.

We can therefore apply Proposition~\ref{Propn std split of trim set}\ref{enu:Prop1.2}
repeatedly to obtain a partition $\widetilde{\mathscr{A}}$ of $C$
such that $\widetilde{\mathscr{A}}\sim\mathscr{B}_{3}\cup\mathscr{B}_{4}=\mathscr{B}_{1}\cup\mathscr{B}_{2}$.
Rearranging, we obtain a partition $\mathscr{A}_{1}\cup\mathscr{A}_{2}$
of $C$ such that elements of $\mathscr{A}_{1}\cup\mathscr{A}_{2}$
are disjoint from each other and $\mathscr{A}_{i}\sim\mathscr{B}_{i}$
for $i=1,2$. Let $C_{i}=\bigcup\mathscr{A}_{i}$, so that $\mu(C_{i})=\mu(D_{i})$
($i=1,2$), to complete this step.

So by Theorem~\ref{(Vaught)} there is an isomorphism $\beta\colon(A)\rightarrow(B)$
such that $\{C\in(A)\mid C\sim C\beta\}$ disjointly generates $(A)$.
Let $\alpha\colon A\rightarrow B$ be the Stone space homeomorphism
induced by $\beta$. If $x\in A\cap X_{p}$, then $x$ has a neighbourhood
base $V_{x}$ of $p$-trim sets by Proposition~\ref{Propn std split of trim set}\ref{enu:Prop1.1},
and we may assume that $C\sim C\beta$ for each $C\in V_{x}$. But
$C$ is $p$-trim iff $\mu(C)=1.p$ iff $C\beta$ is $p$-trim. Hence
$x\alpha$ has a neighbourhood base of $p$-trim sets, so $x\alpha\in B\cap Y_{p}$
(Proposition~\ref{Propn std split of trim set}\ref{enu:Prop1.1}
again), and $(X_{p}\cap A)\alpha\subseteq Y_{p}\cap B$. Equality
follows by considering $\alpha^{-1}$, so $\alpha$ is a $P$-homeomorphism,
as required.
\end{svmultproof}

\section{Distinguished ideals arising from a complete $P$-partition}

We need the following result both for the ``only if'' proof in Theorem~\ref{Thm: P-partition categorical}
below, and for the proof of Theorem~\ref{Theorem: main result}.
\begin{proposition}
\label{Prop:finitely based}Let $R$ be a countable Boolean ring with
Stone space $X$, and $\mathscr{C}$ a closure subalgebra of $2^{X}$
generated by its closed elements such that the system $(R,I(\overline{\mathscr{C}}))$
is $\omega$-categorical. Then $\mathscr{C}$ is finite and each of
its atoms has only finitely many isolated points.
\end{proposition}
\begin{svmultproof}
Let $\varPsi$ be the group of homeomorphisms of $X$ that fix each
closed element of $\mathscr{C}$. Now $\varPsi$ is almost $1$-transitive
on the compact open subsets of $X$, and so there are only finitely
many (say $m$) $\varPsi$-invariant closed subsets of $X$, since
such are complements of unions of $\varPsi$-orbits on the compact
opens. But $\mathscr{C}$ is generated as a Boolean algebra by its
closed elements, and so $\mathscr{C}$ is finite with at most $2^{m}$
elements and $\varPsi$ fixes every element of $\mathscr{C}$. 

If now $A$ is an atom of $\mathscr{C}$ and has at least $k$ isolated
points, then we can find compact open sets $\{B_{j}\mid j\leqslant k\}$
such that $|B_{j}\cap A|=j$, and these belong to different orbits
of $\varPsi$ as $A$ is $\varPsi$-invariant. Hence $k$ is bounded
and $A$ has only finitely many isolated points, as required. 
\end{svmultproof}

To establish the $\omega$-categoricity of many Boolean structures
(in this instance, those arising from a finite-crowded $P$-partition
for some PO system $P$), it suffices to show that the structure-preserving
automorphisms are almost $1$-transitive, which will follow in our
case from Theorem~\ref{Thm:same mu=00003Dhomeom}. This idea forms
the basis of the proof of Theorem~\ref{Thm: P-partition categorical}. 
\begin{theorem}
\label{Thm: P-partition categorical}Let $P$ be a PO system and $\mathscr{X}=\{X_{p}\mid p\in P\}$
a complete $P$-partition of the Stone space $X$ of the countable
Boolean ring $R$. Then the system $(R,\{I(\overline{X_{p}})\mid p\in P\})$
is $\omega$-categorical iff $\mathscr{X}$ is finite-crowded.
\end{theorem}
\begin{remark}
The ``if'' statement for the case of compact $\omega$-Stone spaces
follows from Apps~\cite[Theorem~B]{Apps-Stone}.
\end{remark}
\begin{svmultproof}
Let $L=\{p\in P\mid\overline{X_{p}}\text{ is compact}\}$, let $I_{p}=I(\overline{X_{p}})$,
and let $\mathscr{C}$ be the closure subalgebra of $2^{X}$ generated
by $\{\overline{X_{p}}\mid p\in P\}$.

If the system $(R,\{I_{p}\mid p\in P\})$ is $\omega$-categorical,
then so is the system $(R,\{I(\overline{\mathscr{C}}))$. Hence by
Proposition~\ref{Prop:finitely based} $\mathscr{C}$ is finite and
each of its atoms has only finitely many isolated points. But the
map $P\rightarrow\mathscr{C}:p\mapsto\overline{X_{p}}$ is injective,
since $<$ is antisymmetric on~$P$, and so $P$ is finite. Hence
also $X_{p}=\overline{X_{p}}-\bigcup_{q\lvertneqq p}\overline{X_{q}}\in\mathscr{C}$
and $\mathscr{X}=\At(\mathscr{C})$. Moreover, if $p<p$ then $X_{p}\subseteq X_{p}^{\prime}$
and $X_{p}$ is crowded, while if $p\in P^{d}$ then $X_{p}$ is discrete
and so finite. Therefore $\mathscr{X}$ is finite-crowded.

Now suppose that $\mathscr{X}$ is finite-crowded. We must show that
the automorphisms of $R$ that preserve each $I_{p}$ are almost $r$-transitive
on $R$ for each $r\geqslant1$. Let $X_{L}=\bigcup_{p\in L}\overline{X_{p}}$,
which is closed and compact, and choose $E\in R$ such that $X_{L}\subseteq E$.
Fix $r\geqslant1$ and let $\boldsymbol{A}^{(r)}\in R^{r}$. Set $A_{r+1}=E$
and let $\boldsymbol{B}^{\{s\}}$ be the atoms of the (finite) Boolean
sub-ring of $R$ generated by $\boldsymbol{A}^{\{r+1\}}$, so that
$s\leqslant2^{r+1}$. Let $K_{j}=\{i\leqslant s\mid B_{i}\subseteq A_{j}\}$
for $j\leqslant r+1$, and associate with $\boldsymbol{A}^{(r)}$
the $(r+s+1)$-tuple $\{\mu(B_{1}),\ldots,\mu(B_{s}),\boldsymbol{K}^{(r+1)}\}$.
There are only finitely many possibilities for each $\mu(B_{i})$,
as if $p\in P^{d}$ then $|B_{i}\cap X_{p}|\leqslant|X_{p}|$, and
the latter is finite. Hence there are also only finitely many possibilities
for this $(r+s+1)$-tuple.

Suppose now that $\boldsymbol{C}^{(r)}$ is another $r$-tuple in
$R^{r}$, giving rise to atoms $\boldsymbol{D}^{\{s\}}$ of the Boolean
sub-ring of $R$ generated by $\boldsymbol{C}^{\{r+1\}}$, where $C_{r+1}=E$,
and to the same $(r+s+1)$-tuple, so that $\mu(D_{i})=\mu(B_{i})$
and $B_{i}\subseteq A_{j}$ iff $D_{i}\subseteq C_{j}$ for each $i\leqslant s$
and $j\leqslant r+1$. 

If $X$ is compact, let $u=r+1$ and $v=s$. Then $\bigcup_{j\leqslant u}A_{j}=\bigcup_{j\leqslant u}C_{j}=X$,
as $X_{L}=E=X$. 

If $X$ is not compact, let $u=r+2$. As $X_{L}\subseteq E$ and $X_{p}$
is not compact for $p\notin L$, we can choose $A_{u}$ and $C_{u}$
in $R$, disjoint from $A_{j}$ and $C_{j}$ respectively for $j\leqslant r+1$,
such that $T(A_{u})=T(C_{u})=P-L$. Extending $A_{u}$ and $C_{u}$
as necessary, we may further assume that $\bigcup_{j\leqslant u}A_{j}=\bigcup_{j\leqslant u}C_{j}$.
Let $v=s+1$, $B_{v}=A_{u}$ and $D_{v}=C_{u}$, so that $K_{u}=\{v\}$.
Then $\mu(B_{v})=\mu(D_{v})=\sum_{p\in F}1.p$, where $F=(P-L)_{\min}$,
as $P^{d}\subseteq L$.

For each $i\leqslant v$, apply Theorem~\ref{Thm:same mu=00003Dhomeom}
to find a $P$-homeomorphism $\beta_{i}\colon B_{i}\rightarrow D_{i}$.
These can be combined to give a $P$-homeomorphism $\beta$ of $X$,
with $B_{i}\beta=D_{i}$ for $i\leqslant v$ and $x\beta=x$ for $x\notin\bigcup_{j\leqslant u}A_{j}$.
Let $\alpha$ be the corresponding automorphism of $R$; then $A_{j}\alpha=C_{j}$
for $j\leqslant r$, and $\alpha$ will preserve each $I_{p}$ as
$\beta$ preserves each $\overline{X_{p}}$. Hence the system $(R,\{I_{p}\mid p\in P\})$
is $\omega$-categorical, as required.
\end{svmultproof}

\section{The main result}

We will need the following Lemmas for the proofs both of Theorem~\ref{Theorem: main result}
and of the classification results in Sections~\ref{sec:Classification-POsystem}
and~\ref{sec:Classificationposet}. 
\begin{description}
\item [{\emph{Notation}:}] Let $\mathscr{D}$ be a TBA\@. For $A\in\mathscr{D}$,
let $A^{c}=A\cap A^{\prime}$ and $A^{d}=A-A^{c}$. If $X$ is a topological
space and $\mathscr{D}$ is a sub-TBA of $2^{X}$, then $A^{d}$ is
the set of isolated points in $A$. 
\end{description}
\begin{lemma}
\label{Lemma atoms}Let $\mathscr{D}$ be a TBA and $\mathscr{C}$
a finite closure subalgebra of $(\mathscr{D},^{-})$ such that $\mathscr{D}=\langle\mathscr{C},^{\prime}\rangle$.
Suppose also that $(A^{d})^{\prime}=\emptyset$ for each $A\in\At(\mathscr{C})$.
Then $A^{c}\subseteq(A^{c})^{\prime}=\overline{A}-A^{d}$ for each
$A\in\At(\mathscr{C})$, $\mathscr{D}$ is finite and $\At(\mathscr{D})=\{A^{c},A^{d}\mid A\in\At(\mathscr{C})\}-\{\emptyset\}$. 
\end{lemma}
\begin{svmultproof}
For $A\in\At(\mathscr{C})$, we have $A=A^{c}\dotplus A^{d}$, so
$A^{c}\subseteq A^{\prime}=(A^{c})^{\prime}$, $A^{\prime}\cap A^{d}=A\cap A^{\prime}\cap A^{d}=\emptyset$
and so also $\overline{A^{c}}\cap A^{d}=\emptyset$. Further, $A^{d}$
is closed and $\overline{A}=\overline{A^{c}}\dotplus A^{d}$.

Now let $\mathscr{E}$ be the (finite) Boolean subalgebra of $\mathscr{D}$
generated by the non-empty elements of $\{A^{c},A^{d}\mid A\in\At(\mathscr{C})\}$,
which are all disjoint. Then $\mathscr{E}$ is closed under~$^{\prime}$,
as $(A^{c})^{\prime}=\overline{A^{c}}=\overline{A}-A^{d}$, and $\overline{A}\in\mathscr{C}$
for $A\in\At(\mathscr{C})$. But $\mathscr{C}\subseteq\mathscr{E}$,
and so $\mathscr{D}=\mathscr{E}$, $\mathscr{D}$ is finite, and its
atoms are as required.
\end{svmultproof}

The next Lemma is a consequence of Pierce~\cite[Proposition~8.12]{Pierce}
and Apps~\cite[Lemma~3.3]{AppsBP}; we include a proof which is based
on Naturman~\cite[Lemma~4.2.7]{Naturman}.
\begin{lemma}
\label{Lemma: closure of atoms}Let $\mathscr{C}$ be a finite closure
algebra or TBA which is generated as a Boolean algebra by its closed
elements. If $U$ and $V$ are distinct atoms of $\mathscr{C}$, then
$\overline{U}\neq\overline{V}$.
\end{lemma}
\begin{svmultproof}
A check shows that elements of $\mathscr{C}$ have the form $\bigcup_{1\leqslant i\leqslant r}(C_{i}\cap(X-D_{i}))$,
where $C_{i}$ and $D_{i}$ are closed subsets of $X$ with $C_{i},D_{i}\in\mathscr{C}$,
as $\mathscr{C}$ is generated as a Boolean algebra by its closed
elements. Now $U$ is an atom of $\mathscr{C}$, hence $U\subseteq C\cap(X-D)$
for some closed sets $C,D\in\mathscr{C}$. But $V\cap U=\emptyset$
and so $V\subseteq(X-C)\cup D$. Hence either $V\subseteq X-C$ and
$U\subseteq C$, when $V\cap\overline{U}=\emptyset$, or $V\subseteq D$
and $U\subseteq X-D$, when $U\cap\overline{V}=\emptyset$. In either
case we have $\overline{U}\neq\overline{V}$, as required.
\end{svmultproof}

Our third Lemma illustrates the duality between finite PO systems
and finite TBAs generated by their closed elements. If $P$ is a PO
system, we note that $(2^{P},^{\prime})$ becomes a TBA by setting
$Q^{\prime}=\{p\in P\mid p<q,\text{ some \ensuremath{q\in Q\}}}$
for $Q\subseteq P$. 
\begin{lemma}
\label{Lemma: TBA - PO system equivalence}Let $\mathscr{D}$ be a
finite TBA generated by its closed elements and write $X=1_{\mathscr{D}}$.
Then $\At(\mathscr{D})$ is a complete $P$-partition, $\{X_{p}\mid p\in P\}$
say, of $X$ for some finite PO system $P$. The map $\gamma\colon Q\mapsto\bigcup_{q\in Q}X_{q}$
gives an isomorphism from the TBA $(2^{P},^{\prime})$ to~$\mathscr{D}$,
with lower subsets of $P$ corresponding to closed elements of~$\mathscr{D}$.
Moreover, if $\boldsymbol{Q}^{(n)}$ are lower subsets of $P$, then
$\langle\boldsymbol{Q}^{\{n\}},^{\prime}\rangle=2^{P}$ iff $\langle Q_{1}\gamma,\ldots,Q_{n}\gamma,^{\prime}\rangle=\mathscr{D}$.
\end{lemma}
\begin{svmultproof}
Let $\At(\mathscr{D})=\{X_{p}\mid p\in P\}$, which is a finite partition
of $X$ as $X\in\mathscr{D}$. Define a relation $<$ on $P$ via
$q<p$ iff $X_{q}\subseteq X_{p}^{\prime}$. Now $X_{q}\subseteq\overline{X_{p}}$
if $q<p$. So $<$ is antisymmetric and $(P,<)$ is a PO system by
Lemma~\ref{Lemma: closure of atoms}. Moreover, $X_{p}^{\prime}=\bigcup_{q<p}X_{q}$
as $X_{p}^{\prime}\in\mathscr{D}$, and so $\{X_{p}\mid p\in P\}$
is a complete $P$-partition of $X$. The remaining statements follow
easily, using the fact that $\overline{X_{p}}=\bigcup_{q\leqslant p}X_{q}$.
\end{svmultproof}

We now have all the ingredients for our main result:
\begin{theorem}
\label{Theorem: main result}Let $\boldsymbol{J}^{(n)}$ be distinguished
ideals of the countable Boolean ring $R$, corresponding to closed
subsets $\boldsymbol{C}^{(n)}$ of the Stone space $X$ of $R$. Then
the following are equivalent:
\begin{enumerate}
\item \label{enu:Thm4.1}$(R,\boldsymbol{J}^{(n)})$ is $\omega$-categorical; 
\item \label{enu:Thm4.2}$\langle\boldsymbol{C}^{\{n\}},^{-}\rangle$ is
finite, and each of its atoms has only finitely many isolated points;
\item \label{enu:Thm4.3}$\langle\boldsymbol{C}^{\{n\}},^{\prime}\rangle$
is finite, and its atoms form a finite-crowded partition of $X$;
\item \label{enu:Thm4.4}the atoms of $\langle\boldsymbol{C}^{\{n\}},^{\prime}\rangle$
form a finite-crowded complete $P$-partition of $X$ for some finite
PO system $P$;
\item \label{enu:Thm4.5}$\boldsymbol{C}^{\{n\}}$ can be refined to a finite-crowded
complete $P$-partition of $X$ for some finite PO system $P$.
\end{enumerate}
\end{theorem}
\begin{remark}
If we restrict to compact $\omega$-Stone spaces, when $R$ has a
$1$, the equivalence of~\ref{enu:Thm4.1} and~\ref{enu:Thm4.2}
follows from the original result of Macintyre and Rosenstein~\cite{MacRos},
while the equivalence of statements~\ref{enu:Thm4.1} to~\ref{enu:Thm4.5}
can be deduced from Apps~\cite{AppsBP}.
\end{remark}
\begin{svmultproof}
Let $\mathscr{C}=\langle\boldsymbol{C}^{\{n\}},^{-}\rangle$. 

\ref{enu:Thm4.1}~$\Rightarrow$~\ref{enu:Thm4.2}. Immediate from
Proposition~\ref{Prop:finitely based}.

\ref{enu:Thm4.2}~$\Rightarrow$~\ref{enu:Thm4.3}. Let $\mathscr{D}$
be the sub-TBA of $2^{X}$ generated by $\mathscr{C}$. Then $A^{d}$
is finite for each $A\in\At(\mathscr{C})$. So by Lemma~\ref{Lemma atoms},
$\mathscr{D}$ is finite, and its atoms form a finite-crowded partition
of $X$, as each non-empty $A^{c}$ is crowded for $A\in\At(\mathscr{C})$.

\ref{enu:Thm4.3}~$\Rightarrow$~\ref{enu:Thm4.4}. Follows from
Lemma~\ref{Lemma: TBA - PO system equivalence} 

\ref{enu:Thm4.4}~$\Rightarrow$~\ref{enu:Thm4.5}. Immediate.

\ref{enu:Thm4.5}~$\Rightarrow$~\ref{enu:Thm4.1}. Suppose $\boldsymbol{C}^{\{n\}}$
can be refined to a finite-crowded complete $P$-partition $\{X_{p}\mid p\in P\}$
of $X$ for some finite PO system $P$. Let $I_{p}=I(\overline{X_{p}})$.
Then each $C_{j}$ is a union of some of the $\overline{X_{p}}$,
so any automorphism of $R$ fixing each $I_{p}$ will also fix each~$J_{j}$.
By Theorem~\ref{Thm: P-partition categorical} the system $(R,\{I_{p}\mid p\in P\})$
is $\omega$-categorical, and hence so also is the system $(R,\boldsymbol{J}^{(n)})$. 
\end{svmultproof}

\section{\label{sec:Classification-POsystem}Classification in terms of PO
systems}

Theorem~\ref{Theorem: main result} now enables us to classify $\omega$-categorical
Boolean rings with distinguished ideals in terms firstly of finite
PO systems and secondly of finite posets. We need the following additional
definitions. 
\begin{definition}
If $\boldsymbol{J}^{(n)}$ and $\boldsymbol{K}^{(n)}$ are distinguished
ideals of the Boolean rings $R$ and $S$ respectively, then $(R,\boldsymbol{J}^{(n)})$
and $(S,\boldsymbol{K}^{(n)})$ are \emph{isomorphic }if there is
an isomorphism $\alpha\colon R\rightarrow S$ such that $J_{i}\alpha=K_{i}$
for $i\leqslant n$. Let $\boldsymbol{CB[n]}$ denote the set of isomorphism
classes of countable $\omega$-categorical Boolean rings with $n$
distinguished ideals. 

An \emph{extended PO system }is a triple $(P,L,f)$, where $P$ is
a PO system, $L$ is a lower subset of $P$, and $f\colon L_{\min}^{d}\rightarrow\mathbb{N}_{+}$,
where $L_{\min}^{d}=L_{\min}\cap P^{d}$. We will say that $(P,L,f)$
is a \emph{finite-crowded extended PO system }if $P$ is finite and
$P^{d}\subseteq L_{\min}$.

If $(P,L,f)$ is an extended PO system, then $\mathscr{X}$ is a \emph{complete
$(P,L,f)$-partition of the $\omega$-Stone space $X$ }if it is a
complete $P$-partition of $X$ such that $\overline{X_{p}}$ is compact
iff $p\in L$ and $|X_{p}|=f(p)$ for $p\in L_{\min}^{d}$ (by Remark~\ref{rem:trimpartitions},
$X_{p}$ is finite iff $p\in L_{\min}^{d}$).

A \emph{$Q[n]$-tuple }is a tuple $(P,L,f,\boldsymbol{Q}^{(n)})$
such that $(P,L,f)$ is a finite-crowded extended PO system, each
$Q_{i}$ is a lower subset of $P$, and $\langle\boldsymbol{Q}^{\{n\}},^{\prime}\rangle=2^{P}$.
Two $Q[n]$-tuples $(P_{i},L_{i},f_{i},Q_{i1},\ldots,Q_{in})$ ($i=1,2$)
are \emph{isomorphic }if there is an order-isomorphism $\theta\colon P_{1}\rightarrow P_{2}$
such that $L_{1}\theta=L_{2}$, $Q_{1j}\theta=Q_{2j}$ (all $j$)
and $f_{2}(p\theta)=f_{1}(p)$ for $p\in P_{1}$. Let $\boldsymbol{POQ[n]}$
denote the set of isomorphism classes of $Q[n]$-tuples $(P,L,f,\boldsymbol{Q}^{(n)})$.
\end{definition}
Our classification is underpinned by another result from~\cite{Apps-Stone}: 
\begin{theorem}
\label{Thm: PLF partitions unique}(\cite[Corollary~6.2]{Apps-Stone})
Let $(P,L,f)$ be a finite extended PO system. Then there is an $\omega$-Stone
space $X$ that admits a complete $(P,L,f)$-partition $\mathscr{X}$,
and $(X,\mathscr{X})$ is unique up to $P$-homeomorphisms. 
\end{theorem}
Again, in order to make this paper as self-contained as possible,
we include a fresh proof at the end of this section. Our proof of
existence here, given that $P$ is finite, uses a direct construction
largely due to Pierce~\cite{Pierce}, and is considerably simpler
than the proof in~\cite{Apps-Stone} for the general case of infinite~$P$.

Elements of $\boldsymbol{POQ[n]}$ now become the desired invariants
for our first classification:
\begin{theorem}
\label{Thm: PO system classificn}There is a natural bijection between
$\boldsymbol{CB[n]}$ and \textbf{$\boldsymbol{POQ[n]}$}. 
\end{theorem}
\begin{remark}
If we restrict to Boolean algebras, we obtain a natural bijection
between the isomorphism classes of countable $\omega$-categorical
Boolean algebras with $n$ distinguished ideals, and the isomorphism
classes of tuples $(P,f,\boldsymbol{Q}^{(n)})$ where $P$ is a finite
PO system such that $P^{d}\subseteq P_{\min}$, $f\colon P^{d}\rightarrow\mathbb{N}_{+}$
and $\langle\boldsymbol{Q}^{\{n\}},^{\prime}\rangle=2^{P}$. This
result for Boolean algebras is also a consequence of the work of Touraille~\cite{Touraille}.
\end{remark}
\begin{svmultproof}
Suppose first that $\boldsymbol{J}^{(n)}$ are distinguished ideals
of the countable Boolean ring~$R$, corresponding to closed subsets
$\boldsymbol{C}^{(n)}$ of the Stone space $X$ of $R$, such that
the system $(R,\boldsymbol{J}^{(n)})$ is $\omega$-categorical. Let
$\mathscr{D}$ be the sub-TBA of $2^{X}$ generated by $\boldsymbol{C}^{\{n\}}$,
with atoms $\{X_{p}\mid p\in P\}$, so that $P$ is a finite PO system
and $\{X_{p}\}$ is a finite-crowded complete $P$-partition of $X$
(Theorem~\ref{Theorem: main result}). Let $L=\{p\in P\mid\overline{X_{p}}\text{ is compact}\}$,
which is a lower subset of $P$ containing $P^{d}$, and define $f\colon P^{d}\rightarrow\mathbb{N}_{+}$
by $f(p)=|X_{p}|$. With these definitions, $\{X_{p}\}$ forms a complete
$(P,L,f)$-partition of $X$, with $P^{d}\subseteq L_{\min}$, so
that $(P,L,f)$ is a finite-crowded extended PO system. By Lemma~\ref{Lemma: TBA - PO system equivalence},
each closed set $C_{i}$ corresponds to a lower subset $Q_{i}$, say,
of $P$, and $\boldsymbol{Q}^{\{n\}}$ generates $(2^{P},^{\prime})$
as a TBA\@. Hence $(P,L,f,\boldsymbol{Q}^{(n)})$ is a $Q[n]$-tuple.

Isomorphic $\omega$-categorical systems of the form $(R,\boldsymbol{J}^{(n)})$
will give rise to homeomorphic tuples of the form $(X,\boldsymbol{C}^{(n)})$
where each $C_{j}$ is closed, and hence to isomorphic $Q[n]$-tuples.
Moreover, if two $\omega$-categorical systems of form $(R,\boldsymbol{J}^{(n)})$
give rise to isomorphic $Q[n]$-tuples, then the two systems are isomorphic:
as (relabelling as required to obtain the same underlying $Q[n]$-tuple
$(P,L,f,\boldsymbol{Q}^{(n)})$) the underlying Stone spaces with
their associated $P$-partitions are $P$-homeomorphic by Theorem~\ref{Thm: PLF partitions unique},
and this homeomorphism will also map the associated closed sets to
each other. 

It remains to show that every such $Q[n]$-tuple $(P,L,f,\boldsymbol{Q}^{(n)})$
can arise. But if $(P,L,f)$ is a finite-crowded extended PO system,
then by Theorem~\ref{Thm: PLF partitions unique} we can find an
$\omega$-Stone space~$X$ and a complete $(P,L,f)$-partition $\{X_{p}\mid p\in P\}$
of~$X$. By Lemma~\ref{Lemma: TBA - PO system equivalence}, $\boldsymbol{Q}^{(n)}$
gives rise to closed subsets $\boldsymbol{C}^{(n)}$ of $X$ such
that the atoms of $\langle\boldsymbol{C}^{\{n\}},^{\prime}\rangle$
are $\{X_{p}\}$. In addition, $\{X_{p}\}$ is finite-crowded, as
$P^{d}\subseteq L_{\min}$ and therefore discrete sets are compact
and finite. Applying Theorem~\ref{Theorem: main result}, the system
$(R,\{I(C_{j})\mid j\leqslant n\})$ is $\omega$-categorical and
maps onto $(P,L,f,\boldsymbol{Q}^{(n)})$, as required. 
\end{svmultproof}

\subsection{Theorem~\ref{Thm: PLF partitions unique}: proof of uniqueness}

The original proof of uniqueness in~\cite{Apps-Stone} used a back
and forth argument. Here, we employ a direct decomposition instead.
\begin{svmultproof}
\emph{of uniqueness in Theorem~\ref{Thm: PLF partitions unique}}

Suppose that $(P,L,f)$ is a finite extended PO system and that $\{X_{p}\mid p\in P\}$
is a complete $(P,L,f)$-partition of the Stone space $X$ of the
countable Boolean ring $R$. Let $F=L_{\min}$, $G=(P-L)_{\min}$
and $Q=F\cup G$. For uniqueness, it suffices to show that we can
write $R=\bigoplus{}_{p\in Q}\{\bigoplus_{j<n_{p}+1}(A_{pj})\}$,
where $A_{pj}$ is $p$-trim and $n_{p}=1$ if $p\in L_{\min}-P^{d}$,
$n_{p}=f(p)$ if $p\in L_{\min}^{d}$, and $n_{p}=\infty$ if $p\in G$:
for then $F$, $G$, $Q$ and each $n_{p}$ are determined by $(P,L,f)$,
and any two $p$-trim sets are $P$-homeomorphic (Theorem~\ref{Thm:same mu=00003Dhomeom}).

To see this, find $A\in R$ such that $\bigcup_{p\in L}\overline{X_{p}}\subseteq A$.
Using Proposition~\ref{Propn std split of trim set}\ref{enu:Prop1.4},
write $A=A_{1}\dotplus\ldots\dotplus A_{n}$, where each $A_{j}$
is trim. Discarding any redundant $A_{j}$'s as necessary and reducing
$A$ accordingly, we may assume that $T(A_{j})\cap L\neq\emptyset$
for each $j$, so that $L\subseteq T(A)$ and $T(A)_{\min}=L_{\min}=F$.
Then $A$ admits a $\mu(A)$-partition, where $\mu(A)=\sum_{p\in F}n_{p}.p$
and $n_{p}$ are as above, as $|A\cap X_{p}|=f(p)$ for $p\in F\cap P^{d}=L_{\min}^{d}$.

If $L=P$, then $A=X$ and we are done. Otherwise, let $S=\{B\in R\mid B\cap A=\emptyset\}$,
so that $R=(A)\oplus S$. By a routine argument we can find disjoint
$\{B_{n}\in S\mid n\geqslant1\}$ that generate $S$, with $T(B_{n})\subseteq P-L$
for all $n$. As $\overline{X_{p}}$ is not compact for $p\in P-L$,
we can find $m_{1}$ such that $T(B_{1}\cup\cdots\cup B_{m_{1}})=P-L$.
By ``clumping together'' the $B_{n}$ in this way, we can find disjoint
$\{C_{n}\in S\mid n\geqslant1\}$ that generate $S$ such that $T(C_{n})=P-L$
for each $n$. Now split each $C_{n}$ into $q$-trim sets, one for
each $q\in G$, to give the desired decomposition of $S$ as $\bigoplus_{p\in G}\{\bigoplus_{j\geqslant1}(A_{pj})\}$,
where each $A_{pj}$ is $p$-trim.
\end{svmultproof}

\subsection{Theorem~\ref{Thm: PLF partitions unique}: proof of existence}

Let $D_{1}$ and $D_{0}$ denote the Cantor set and the Cantor set
minus a point respectively, whose associated Boolean rings are countable
atomless with and without a $1$ respectively. The following Lemma
and proof of existence of complete $(P,L,f)$-partitions for finite
$P$ are based on proofs by Pierce~\cite[Theorem~4.6]{Pierce} which
handle the compact cases; we have extended these to cover the non-compact
cases as well.
\begin{lemma}
\label{Lemma: Pierce}Let $W$ be an $\omega$-Stone space and $C$
a non-empty closed subset of $W$. Then we can find an $\omega$-Stone
space $X$ and a homeomorphism $\alpha\colon W\rightarrow W\alpha\subseteq X$
such that $W\alpha$ is closed in $X$ and $\overline{Y}\cap W\alpha=C\alpha$
where $Y=X-W\alpha$. Moreover, we can arrange for $Y$ to be infinite
discrete or crowded, and for $\overline{Y}$ to be either compact
(provided $C$ is compact) or non-compact.
\end{lemma}
\begin{svmultproof}
\textbf{Case 1}: $Y$ is to be crowded.

Let $D=D_{0}$ if $C$ is compact and $\overline{Y}$ is to be compact
(case 1A), and let $D=D_{1}$ if $\overline{Y}$ is to be non-compact
(case 1B). Select $x_{0}\in D$, and define the subset $X$ of $W\times D$
by $X=(W\times\{x_{0}\})\cup(C\times D)$; this is a closed subset
of a Stone space, and so is itself a Stone space. Define $\alpha\colon W\rightarrow X:y\mapsto(y,x_{0})$
which is a homeomorphism onto $W\times\{x_{0}\}$. We have $Y=X-W\alpha=C\times(D-\{x_{0}\})$,
so $\overline{Y}=C\times D$, which is compact (case 1A) or non-compact
(case 1B), $\overline{Y}\cap W\alpha=C\times\{x_{0}\}=C\alpha$, and
neither $D-\{x_{0}\}$ nor $Y$ has any isolated points, as required.

\textbf{Case 2}: $Y$ is to be infinite discrete.

Let $V=\emptyset$ if $C$ is compact and $\overline{Y}$ is to be
compact (case 2A), and let $V=\{2,3,\ldots\}$ if $\overline{Y}$
is to be non-compact (case 2B). Let $\{x_{n}\mid n\geqslant1\}$ be
a countable dense discrete subset of $C$, let $Z=\{0\}\cup\{\frac{1}{n}\mid n\geqslant1\}\cup V$,
a Stone space which is compact in case 2A, and define the closed subset
$X$ of $W\times Z$ by $X=(W\times\{0\})\cup Y$, where

\[
Y=\{(x_{m},\frac{1}{n})\mid m<n\}\cup(\{x_{1}\}\times V)
\]

Define $\alpha\colon W\rightarrow X:y\mapsto(y,0)$ which is a homeomorphism
onto $W\times\{0\}$, with $Y=X-W\alpha$ being countably discrete.
Then $\overline{Y}=(C\times\{0\})\cup Y$, $\overline{Y}\cap W\alpha=C\alpha$,
and $\overline{Y}$ is compact in case 2A (being a closed subset of
the compact $C\times Z$) and is non compact in case 2B.

Finally, each of the spaces $X$ so constructed is a closed subset
of a product of Stone spaces, each with a countable base, and so itself
has a countable base; hence each $X$ is an $\omega$-Stone space.
\end{svmultproof}

This Lemma enables the main induction step required to construct an
$\omega$-Stone space that admits a complete $(P,L,f)$-partition.
\begin{svmultproof}
\emph{of existence in Theorem~\ref{Thm: PLF partitions unique}}:
by induction on $|P|$. 

If $|P|=1$, there are 4 cases. Let $P=\{p\}$.

If $P^{d}=\emptyset$, we can take $X$ to be $D_{1}$ if $L=P$ or
$D_{0}$ if $L=\emptyset$. 

If $P^{d}=\{p\}$, we can take $X$ to be a set of $f(p)$ discrete
points if $L=P$ or a countably infinite set of discrete points if
$L=\emptyset$.

Suppose now that the result holds for $|P|\leqslant k$, and that
$|P|=k+1$. Pick a maximal $p\in P$; let $Q=P-\{p\}$, $M=L-\{p\}$,
and let $g$ be the restriction of $f$ to $M_{\min}\cap Q^{d}\subseteq L_{\min}^{d}$.
By the induction hypothesis, we can find an $\omega$-Stone space
$W$ that admits a complete $(Q,M,g)$-partition $\{Y_{q}\mid q\in Q\}$.

Let $Q_{1}=\{q\in Q\mid q<p\}$ and $C=\bigcup_{q\in Q_{1}}Y_{q}$,
which is closed in $W$. If $Q_{1}=\emptyset$ (i.e. $p\in P_{\min}$),
we can take $X=W\dotplus Y_{p}$ with partition $\{Y_{q}\mid q\in P\}$,
where $Y_{p}$ is determined as for when $P=\{p\}$.

Otherwise, use Lemma~\ref{Lemma: Pierce} to find an $\omega$-Stone
space $X$ and a homeomorphism $\alpha\colon W\rightarrow W\alpha\subseteq X$
such that

$\begin{cases}
W\alpha\text{ is closed in }X\\
\overline{X_{p}}\cap W\alpha=C\alpha\text{, where }X_{p}=X-W\alpha\\
\overline{X_{p}}\text{ is compact iff }p\in L\text{ (as if }p\in L\text{, then }Q_{1}\subseteq L\text{ and so }C\text{ is compact)}\\
X_{p}\text{ is crowded if }p\notin P^{d}\text{ and }X_{p}\text{ is infinite and discrete if }p\in P^{d}.
\end{cases}$

Let $X_{q}=Y_{q}\alpha$ for $q\neq p$. An easy check now shows that
$\{X_{q}\mid q\in P\}$ is the required complete $(P,L,f)$-partition
of $X$.
\end{svmultproof}

\section{\label{sec:Classificationposet}Classification in terms of posets}

The $Q[n]$-tuples for our classification in terms of PO systems captured
four aspects of the atoms $\mathscr{X}$ of $\langle\boldsymbol{C}^{\{n\}},^{\prime}\rangle$
(using the notation of Theorem~\ref{Theorem: main result}): namely,
the derived set behaviour of $\mathscr{X}$; which elements of $\mathscr{X}$
are relatively compact; the size of any finite elements of $\mathscr{X}$;
and the mapping from $\boldsymbol{C}^{(n)}$ into $2^{\mathscr{X}}$.

For the classification in terms of posets we need to capture five
aspects of the atoms $\mathscr{Y}$ of $\langle\boldsymbol{C}^{\{n\}},^{-}\rangle$:
namely, the closure behaviour of $\mathscr{Y}$; which elements of
$\mathscr{Y}$ are relatively compact; which elements of $\mathscr{Y}$
are finite; the number of isolated points in each element of $\mathscr{Y}$;
and the mapping from $\boldsymbol{C}^{(n)}$ into $2^{\mathscr{Y}}$.
\begin{definition}
A \emph{$U[n]$-tuple }is a tuple $(S,M,F,g,\boldsymbol{U}^{(n)})$
where $S$ is a finite poset and such that:
\begin{enumerate}
\item $M,U_{1},\ldots,U_{n}$ are lower subsets of~$S$;
\item $F\subseteq M_{\min}$;
\item $g\colon S\rightarrow\mathbb{N}$ is such that $g(s)>0$ for $s\in F$;
\item $\langle\boldsymbol{U}^{\{n\}},^{-}\rangle=2^{S}$. 
\end{enumerate}
We say that two $U[n]$-tuples $(S_{i},M_{i},F_{i},g_{i},U_{i1},\ldots,U_{in})$
($i=1,2$) are \emph{isomorphic }if there is an order-isomorphism
$\theta\colon S_{1}\rightarrow S_{2}$ such that $M_{1}\theta=M_{2}$,
$F_{1}\theta=F_{2}$, $U_{1j}\theta=U_{2j}$ (all $j$) and $g_{2}(s\theta)=g_{1}(s)$
for $s\in S_{1}$. Let $\boldsymbol{PU[n]}$ denote the set of isomorphism
classes of $U[n]$-tuples.
\end{definition}
\begin{theorem}
\label{Thm: classificn posets}There is a natural bijection between
$\boldsymbol{CB[n]}$ and $\boldsymbol{PU[n]}$.
\end{theorem}
\begin{remark}
It follows that there is a natural bijection between elements of \textbf{$\boldsymbol{POQ[n]}$
}and of $\boldsymbol{PU[n]}$. This can also be shown directly using
the approach of Lemma~\ref{Lemma: split poset into PO system}.
\begin{remark}
\label{rem:Alaev}We note that for $D\in\At(\mathscr{C})$, where
$\mathscr{C}$ is the (finite) closure subalgebra of $2^{X}$ generated
by $\boldsymbol{C}^{\{n\}}$, $D$ and $\overline{D}$ have the same
number of isolated points and are both either finite, relatively compact
or not relatively compact. Using the duality between finite Heyting
algebras and finite closure algebras generated by their closed elements,
we recover the result of Alaev~\cite[Theorem~2.5]{Alaev}, that a
countable $\omega$-categorical Boolean algebra $R$ with distinguished
ideals $\boldsymbol{J}^{(n)}$ is uniquely determined by the isomorphism
type of $\widehat{\mathscr{C}}$, the Heyting sub-algebra of $H(R)$
generated by $\boldsymbol{J}^{\{n\}}$ and with distinguished elements
$\boldsymbol{J}^{(n)}$, by the number of atoms in $R/J$ for each
$J\in\widehat{\mathscr{C}}$, and whether each such $R/J$ is finite
or infinite.
\begin{remark}
Macintyre and Rosenstein~\cite{MacRos} investigated countable $\omega$-categorical
``totally atomless augmented Boolean algebras'' of the form $(R,\mathscr{H})$,
where $\mathscr{H}$ is a finite Heyting subalgebra of $H(R)$ such
that $R/J$ is atomless for each $J\in\mathscr{H}$, and each $J\in\mathscr{H}$
is distinguished. These correspond exactly to isomorphism classes
of $U[n]$-tuples $(S,M,F,g,\boldsymbol{U}^{(n)})$ where $S$ is
the finite poset dually equivalent to $\mathscr{H}$, $M=2^{S}$,
$F=\emptyset$, $g(s)=0$ for all $s\in S$, and $\boldsymbol{U}^{\{n\}}$
comprises all the lower subsets of $S$. So if we factor out the order
in which the distinguished elements of $\mathscr{H}$ appear, we obtain
a natural bijection between the isomorphism classes of countable $\omega$-categorical
totally atomless augmented Boolean algebras and the isomorphism classes
of finite posets (equivalently, of finite PO systems $P$ such that
$P^{d}=\emptyset$).
\end{remark}
\end{remark}
\end{remark}
\begin{svmultproof}
\textbf{Step 1: define a map $\hat{\phi}$ from $\boldsymbol{CB[n]}$
into $\boldsymbol{PU[n]}$}

Let $\boldsymbol{J}^{(n)}$ be distinguished ideals of the countable
Boolean ring $R$, corresponding to closed subsets $\boldsymbol{C}^{(n)}$
of the Stone space $X$, such that $(R,\boldsymbol{J}^{(n)})$ is
$\omega$-categorical. Let $\mathscr{C}$ be the (finite) closure
subalgebra of $2^{X}$ generated by $\boldsymbol{C}^{\{n\}}$. Let
$(R,\boldsymbol{J}^{(n)})\phi$ be the tuple $(S,M,F,g,\boldsymbol{U}^{(n)})$,
where:

$\begin{cases}
S=\At(\mathscr{C})\text{, with order relation }A\leqslant B\text{ in }S\text{ if }A\subseteq\overline{B}\text{ in }X;\\
M=\{A\in S\mid\overline{A}\text{ is compact}\};\\
F=\{A\in S\mid A\text{ is finite}\};\\
g(A)\text{ is the number of isolated points in }A\text{, which is finite by Theorem\,\ref{Theorem: main result};}\\
U_{j}=\{A\in S\mid A\subseteq C_{j}\}.
\end{cases}$

Now $S$ is a poset by Lemma~\ref{Lemma: closure of atoms}, and
$\boldsymbol{U}^{\{n\}}$ generates $2^{S}$ since $\mathscr{C}\cong2^{S}$
as closure algebras. Moreover, $M$ and $U_{j}$ are lower subsets
of $S$ as $C_{j}$ is closed, and $F\subseteq M_{\min}$ as $A$
is closed for $A\in F$. Hence $(S,M,F,g,\boldsymbol{U}^{(n)})$ is
a $U[n]$-tuple. It follows easily that isomorphic $\omega$-categorical
systems of the form $(R,\boldsymbol{J}^{(n)})$ will give rise to
isomorphic $U[n]$-tuples in $\boldsymbol{PU[n]}$, so we obtain an
induced map $\hat{\phi}\colon\boldsymbol{CB[n]}\rightarrow\boldsymbol{PU[n]}$.

We must show firstly that an $\omega$-categorical $(R,\boldsymbol{J}^{(n)})$
is determined up to isomorphism by the isomorphism class of the associated
$U[n]$-tuple, and secondly that every possible $U[n]$-tuple in $\boldsymbol{PU[n]}$
can arise. To do this, we will link $U[n]$-tuples with $Q[n]$-tuples
and use the corresponding results for $Q[n]$-tuples.

\textbf{Step 2: $\hat{\phi}$ is injective}

Let $(R,\boldsymbol{J}^{(n)})$, $\boldsymbol{C}^{(n)}$, $\mathscr{C}$
and $(S,M,F,g,\boldsymbol{U}^{(n)})$ be as above. Let $\mathscr{D}=\langle\boldsymbol{C}^{\{n\}},^{\prime}\rangle=\langle\mathscr{C},^{\prime}\rangle$,
a finite sub-TBA of $2^{X}$, and let $P=\At(\mathscr{D})$, viewed
as a PO system with $A<B$ in $P$ iff $A\subseteq B^{\prime}$ in
$X$. Let $(P,L,f,\boldsymbol{Q}^{(n)})$ be the corresponding $Q[n]$-tuple
associated with $(R,\boldsymbol{J}^{(n)})$, as per the construction
in Theorem~\ref{Thm: PO system classificn}. Now $\mathscr{C}$ satisfies
the conditions of Lemma~\ref{Lemma atoms}, and so $\At(\mathscr{D})=\{A^{c},A^{d}\mid A\in\At(\mathscr{C})\}-\{\emptyset\}$.
But $A^{d}$ consists of the (finite number of) isolated points in
$A$ for $A\in\At(\mathscr{C})$, so $A^{d}=\emptyset$ iff $g(A)=0$,
and $A^{c}=\emptyset$ iff $A=A^{d}$ iff $A\in F$. Also, $(A^{c})^{\prime}=\overline{A}-A^{d}$
and $(A^{d})^{\prime}=\emptyset$ for $A\in\At(\mathscr{C})$. So
the $U[n]$-tuple $(S,M,F,g,\boldsymbol{U}^{(n)})$ determines the
isomorphism class of the PO system $(P,<)$, with $P^{d}=\{A^{d}\mid A\in\At(\mathscr{C})\wedge g(A)>0\}$.

For $A\in\At(\mathscr{C})$, we have $f(A^{d})=|A^{d}|=g(A)$ if $g(A)>0$,
and $\overline{A}$ is compact iff $\overline{A^{c}}$ is compact
(as $\overline{A^{c}}=(A^{c})^{\prime}=\overline{A}-A^{d}$). Hence
$L=\{A^{d}\mid g(A)>0\}\cup\{A^{c}\mid A\in M-F\}$. Finally, $C_{j}=\bigcup\{A\mid A\in U_{j}\}$,
so $Q_{j}=\{A^{d}\mid A\in U_{j}\wedge g(A)>0\}\cup\{A^{c}\mid A\in U_{j}-F\}$.

Given $(R,\boldsymbol{J}^{(n)})$, we have shown that the associated
$U[n]$-tuple $(S,M,F,g,\boldsymbol{U}^{(n)})$ determines the isomorphism
class of the associated $Q[n]$-tuple $(P,L,f,\boldsymbol{Q}^{(n)})$.
Hence if two $\omega$-categorical systems of the form $(R,\boldsymbol{J}^{(n)})$
generate isomorphic $U[n]$-tuples, then they generate isomorphic
$Q[n]$-tuples and so are isomorphic by Theorem~\ref{Thm: PO system classificn}.

\textbf{Step 3: $\hat{\phi}$ is surjective }

Suppose now that $(S,M,F,g,\boldsymbol{U}^{(n)})$ is a $U[n]$-tuple.
Let $H=\{s\in S-F\mid g(s)>0\}$. Now $F\subseteq S_{\min}$, so we
can apply Lemma~\ref{Lemma: split poset into PO system} below to
the partition $(F,S-F-H,H)$ of $S$ to obtain a finite PO system
$P$ and a closure algebra morphism $\alpha\colon(2^{S},^{-})\rightarrow(2^{P},^{-})$
satisfying the properties of that Lemma. In particular, in the notation
of Lemma~\ref{Lemma: split poset into PO system} and Lemma~\ref{Lemma atoms},
$P^{d}\cap(\{s\}\alpha)=(\{s\}\alpha)^{d}=\{s_{d}\}$ for $s\in F\cup H$,
and $P^{d}=\{s_{d}\mid s\in F\cup H\}$.

Set $L=P^{d}\cup M\alpha$, which is a lower subset of $P$ since
$M$ is a lower subset of $S$, with $P^{d}\subseteq L_{\min}$. For
$s\in F\cup H$, set $f(s_{d})=g(s)$, and let $Q_{j}=U_{j}\alpha$,
which is a lower subset of $P$ as $U_{j}$ is a lower subset of $S$.
Moreover, $\langle\boldsymbol{Q}^{\{n\}},^{\prime}\rangle=2^{P}$
by Lemma~\ref{Lemma: split poset into PO system}, and so $(P,L,f,\boldsymbol{Q}^{(n)})$
is a $Q[n]$-tuple.

By Theorem~\ref{Thm: PO system classificn}, we can find an $\omega$-categorical
system $(R,\boldsymbol{J}^{(n)})$ that gives rise to $(P,L,f,\boldsymbol{Q}^{(n)})$.
Let $X$ be the Stone space of $R$, let $\boldsymbol{C}^{(n)}$ be
the $n$-tuple of closed subsets of $X$ corresponding to $\boldsymbol{J}^{(n)}$,
and let $\mathscr{D}$ and $\mathscr{C}$ be the sub-TBA and closure
subalgebra of $2^{X}$ respectively generated by $\boldsymbol{C}^{\{n\}}$.
Let $\beta\colon2^{P}\rightarrow\mathscr{D}$ be the corresponding
TBA isomorphism (Lemma~\ref{Lemma: TBA - PO system equivalence},
as $P$ indexes $\At(\mathscr{D})$), so that $|\{s_{d}\}\beta|=f(s_{d})$
for $s\in F\cup H$; for $T\subseteq P$, $\overline{T\beta}$ is
compact iff $T\subseteq L$; and $Q_{j}\beta=C_{j}$. Let $\gamma=\alpha\beta\colon2^{S}\rightarrow\mathscr{D}$,
which is a closure algebra morphism. Then $U_{j}\gamma=Q_{j}\beta=C_{j}$;
hence $2^{S}\gamma=\langle\boldsymbol{U}^{\{n\}},^{-}\rangle\gamma=\langle\boldsymbol{C}^{\{n\}},^{-}\rangle=\mathscr{C}$.
Let $Y_{s}=\{s\}\gamma$ for $s\in S$, so that $\At(\mathscr{C})=\{Y_{s}\mid s\in S\}$. 

Now for $s\in S$, $\overline{Y_{s}}$ is compact iff $\{s\}\alpha\subseteq L$
iff $s\in F\cup M=M$, and $Y_{s}$ is finite iff $\{s\}\alpha\subseteq P^{d}$
iff $s\in F$. Furthermore, for $s\in S$ the isolated points of $Y_{s}$
are $Y_{s}^{d}=[(\{s\}\alpha)^{d}]\beta$, as $\beta$ is a TBA-isomorphism;
so if $s\in F\cup H$ then $Y_{s}^{d}=\{s_{d}\}\beta$ and has size
$f(s_{d})$, while $Y_{s}^{d}=\emptyset$ if $s\notin(F\cup H)$;
so in either case $Y_{s}$ has $g(s)$ isolated points. Hence the
map $S\rightarrow\At(\mathscr{C}):s\mapsto Y_{s}$ is a poset isomorphism,
and it induces an isomorphism mapping the $U[n]$-tuple $(S,M,F,g,\boldsymbol{U}^{(n)})$
to $(R,\boldsymbol{J}^{(n)})\phi$, as required.
\end{svmultproof}

To complete the proof of Theorem~\ref{Thm: classificn posets}, we
need the following Lemma which ``splits the atoms'' of a poset $S$
to obtain a PO system $P$, given a requirement as to which atoms
of $S$ become discrete in $P$, which become crowded in $P$, and
which split into a join of discrete and crowded atoms of $P$. We
use the notation of Lemma~\ref{Lemma atoms}, which this Lemma effectively
complements. As usual, if $P$ is a PO system, then $(2^{P},^{\prime})$
becomes a TBA by setting $Q^{\prime}=\{p\in P\mid p<q,\text{ some \ensuremath{q\in Q\}}}$
and becomes a closure algebra with $\overline{Q}=\{p\in P\mid p\leqslant q,\text{ some \ensuremath{q\in Q\}}}$,
for $Q\subseteq P$.
\begin{lemma}
\label{Lemma: split poset into PO system}Let $S$ be a finite poset,
and $(F,G,H)$ a partition of $S$ such that $F\subseteq S_{\min}$,
where some of $F,G,H$ may be empty. Then there is a finite PO system
$(P,<)$ such that $P^{d}\subseteq P_{\min}$, and a closure algebra
morphism $\alpha\colon(2^{S},^{-})\rightarrow(2^{P},^{-})$ such that
for $s\in S$:
\begin{enumerate}
\item \label{enu:-3.1}$(\{s\}\alpha)^{d}\subseteq P^{d}$;
\item $|(\{s\}\alpha)^{d}|=1$ for $s\in F\cup H$ and $0$ for $s\in G$;
\item $(\{s\}\alpha)^{c}=\emptyset$ iff $s\in F$;
\item $P=\bigcup_{s\in S}\{s\}\alpha$;
\item \label{enu:-3.5}$\langle2^{S}\alpha,^{\prime}\rangle=2^{P}$.
\end{enumerate}
Moreover, if $\{T_{1},\ldots,T_{n}\}$ are subsets of $S$ that generate
$2^{S}$ as a closure algebra, then $\{T_{1}\alpha,\ldots,T_{n}\alpha\}$
generate $2^{P}$ as a TBA\@.
\end{lemma}
\begin{remark}
It follows that $[(\{s\}\alpha)^{d}]^{\prime}=\emptyset$ for $s\in S$,
and so by Lemma~\ref{Lemma atoms} the pair $(P,\alpha)$ is unique
up to isomorphism: that is, if $(Q,\beta)$ is another such pair,
then there is a PO system isomorphism $\theta\colon P\rightarrow Q$
such that $\{s\}\alpha\tilde{\theta}=\{s\}\beta$ for all $s\in S$,
where $\tilde{\theta}$ is the induced isomorphism from $2^{P}$ to
$2^{Q}$.
\end{remark}
\begin{svmultproof}
Let $P=\{s_{d}\mid s\in F\cup H\}\cup\{s_{c}\mid s\in G\cup H\}$,
with order relation for $s\in S$ given by the following (where $p<q$
iff $p\in\{q\}^{\prime}$):

$\begin{cases}
\{s_{d}\}{}^{\prime}=\emptyset\text{ for }s\notin G;\\
\{s_{c}\}^{\prime}=\{r_{c}\mid r\notin F\wedge r\leqslant s\}\cup\{r_{d}\mid r\notin G\wedge r\lneqq s\}\text{ for }s\notin F.
\end{cases}$

Let $\{s\}\alpha=\{s_{d}\}$ or $\{s_{c}\}$ or $\{s_{d},s_{c}\}$
according as $s\in F$ or $s\in G$ or $s\in H$. We note that $(\{s\}\alpha)^{c}=(\{s\}\alpha)\cap(\{s\}\alpha)^{\prime}=\{s_{c}\}$
for $s\in G\cup H$; further, $(\{s\}\alpha)^{d}=\{s_{d}\}$ for $s\in F\cup H$
and is empty for $s\in G$. Hence $\langle2^{S}\alpha,^{\prime}\rangle=2^{P}$
and $P^{d}=\{s_{d}\mid s\in F\cup H\}\subseteq P_{\min}$. We have
now established properties~\ref{enu:-3.1} to~\ref{enu:-3.5}.

With this relation, $(P,<)$ becomes a PO system (transitivity follows
from that for~$S$), and $(2^{P},^{-})$ becomes a closure algebra.
To show that $\alpha$ is a morphism of closure algebras, it suffices
to show that $\overline{\{s\}\alpha}=\overline{\{s\}}\alpha$ in $2^{P}$
for $s\in S$ since $P=\bigcup_{s\in S}\{s\}\alpha$. This follows
easily from the fact that $\overline{\{s\}\alpha}=\{r_{c}\mid r\notin F\wedge r\leqslant s\}\cup\{r_{d}\mid r\notin G\wedge r\leqslant s\}$;
we note that this holds for $s\in F$ as $F\subseteq S_{\min}$. 

Finally, if $\{T_{1},\ldots,T_{n}\}$ are subsets of $S$ that generate
$2^{S}$ as a closure algebra, then $\langle\{T_{1}\alpha,\ldots,T_{n}\alpha\},^{\prime}\rangle$
contains $2^{S}\alpha$, and so equals $2^{P}$, as required.
\end{svmultproof}

\section{Links with previous results\label{sec:Links-with-previous}}

We have used some results from Apps~\cite{Apps-Stone}, namely Theorems~\ref{Thm:same mu=00003Dhomeom}
and~\ref{Thm: PLF partitions unique}, concerning complete finite
$P$-partitions of locally compact Stone spaces of Boolean \uline{rings}
which may or may not have a $1$, and included simplified proofs appropriate
for the case of finite $P$.

The other relevant work to date in this area has concerned Boolean
\uline{algebras}: the case where $X$ is the compact Stone space
of the Boolean algebra $R$, and there are closed subsets $\boldsymbol{C}^{(n)}$
of $X$ corresponding to distinguished ideals $\boldsymbol{J}^{(n)}$
of~$R$. Let $\mathscr{C}$ be the closure subalgebra of $2^{X}$
generated by $\boldsymbol{C}^{\{n\}}$, and let $\widehat{\mathscr{C}}$
be the corresponding Heyting sub-algebra $\{I(\overline{C})\mid C\in\mathscr{C}\}$
of $H(R)$ with distinguished elements $\boldsymbol{J}^{(n)}$. Here,
if $I$ and $J$ are ideals of $R$, then we set $I.J=I\cap J$, $I+J=\{A+B\mid A\in I,B\in J\}$
and $I\rightarrow J=\{A\in R\mid(B\leqslant A)\wedge(B\in I\Rightarrow B\in J)\}$.

The result of Macintyre and Rosenstein~\cite[Theorem~7]{MacRos},
proved in the context of Heyting algebras, translates to the equivalence
of statements~\ref{enu:Thm4.1}~and~\ref{enu:Thm4.2} in Theorem~\ref{Theorem: main result},
as elements of $\widehat{\mathscr{C}}$ correspond exactly to the
closed elements of $\mathscr{C}$; $\widehat{\mathscr{C}}$ is finite
iff $\mathscr{C}$ is finite; and an atom and its closure have the
same isolated points.

Apps~\cite[Theorem~C]{AppsBP}, as part of an investigation into
$\omega$-categorical finite extensions of certain Boolean powers
of groups, showed that that for any finite PO system $P$ there is
a compact $\omega$-Stone space $X$ that admits a complete $P$-partition
$\mathscr{X}=\{X_{p}\mid p\in P\}$ such that $|X_{p}|=1$ if $p\in P^{d}$;
that $\{X,\mathscr{X}\}$ is unique up to $P$-homeomorphisms; and
that the corresponding Boolean algebra system is $\omega$-categorical.
(We note that~\cite[Theorem~C]{AppsBP} concerns the more general
situation of Stone spaces acted on by a finite group; the above stated
result follows by taking the group action to be trivial.) By splitting
each finite partition element into singletons, this result, together
with~\cite[Theorem~B]{AppsBP}, leads to the ``if'' statement in
Theorem~\ref{Thm: P-partition categorical} and to the equivalence
of statements~\ref{enu:Thm4.1}~to~\ref{enu:Thm4.5} of Theorem~\ref{Theorem: main result}
for the case of Boolean algebras. 

A Boolean algebra $R$ is \emph{pseudo-indecomposable (PI)} if for
all $A\in R$, either $(A)\cong R$ or $(1-A)\cong R$; and $B\in R$
is PI if $(B)$ is PI\@. Pal'chunov~\cite{Palchunov} extended this
definition to cover Boolean algebras with distinguished ideals, termed
\emph{$I$-algebras}. He showed that any $\omega$-categorical $I$-algebra
can be decomposed as a finite direct product of PI $I$-algebras,
and he identified propositions which serve as axioms for the PI $I$-algebras. 

Touraille~\cite{Touraille} defined a \emph{Heyting-sa algebra} to
be a Heyting algebra with an additional unary operation $\sa()$ satisfying
certain conditions. For such an algebra $\mathscr{F}$, he defined
$M(\mathscr{F})$ to be the set of maximal elements $J\in\mathscr{F}-\{1\}$
such that $\sa(J)=J$. The Heyting algebra $H(R)$ of ideals of $R$
becomes a Heyting sa-algebra with the definition $\sa(J)=\{B\in R\mid B/J\text{ is atomless in }R/J\}$.
Let $\mathscr{D}_{R}$ be its sub-algebra generated by $\boldsymbol{J}^{\{n\}}$
with distinguished elements $\boldsymbol{J}^{(n)}$, and for $J\in M(\mathscr{D}_{R})$
let $\eta_{R}(J)\in\mathbb{N}_{+}\cup\{\infty\}$ be the number of
atoms in $R/J$. A conseqence of Touraille's work is that the isomorphism
type of an $\omega$-categorical $I$-algebra is uniquely determined
by the pair $(\mathscr{D}_{R},\eta_{R})$, and that every possible
pair $(\mathscr{F},\eta)$ can arise, where $\mathscr{F}$ is a finite
Heyting sa-algebra generated by some distinguished elements and satisfying
an extra condition, and $\eta$ is a function from $M(\mathscr{F})$
to $\mathbb{N}_{+}$ (see~\cite[Theorem~2.2,~ Proposition~2.3]{Alaev}).
Touraille~\cite{Touraille,Touraille2} also showed that finite Heyting
sa-algebras correspond exactly firstly to finite TBAs generated by
their closed elements and secondly to finite PO systems. Under this
correspondence $\mathscr{D}_{R}$ becomes $\mathscr{D}$, the sub-TBA
of $2^{X}$ generated by $\boldsymbol{C}^{\{n\}}$ and with distinguished
elements $\boldsymbol{C}^{(n)}$, corresponding to a finite PO system
$P$; $M(\mathscr{D}_{R})$ becomes the finite minimal closed sets
in~$\mathscr{D}$, corresponding to $P_{\min}^{d}$; and $\eta$
becomes $f\colon P_{\min}^{d}\rightarrow\mathbb{N}_{+}$. Hence Touraille's
result for Heyting sa-algebras yields Theorem~\ref{Thm: PO system classificn}
for the case of Boolean algebras.

Alaev~\cite{Alaev} further developed the results of Pal'chunov and
Touraille and described $\omega$-categorical $I$-algebras via their
decomposition into PI $I$-algebras. He showed that $\omega$-categorical
PI $I$-algebras correspond precisely to pairs of the form $(H,\varepsilon)$,
where $H$ is a finite Heyting algebra with $n$ distinguished non-zero
generating elements (for some $n$), $H-\{1\}$ has a greatest element
$H_{0}$ (say), and $\varepsilon\in\{0,1\}$. Under this correspondence,
the Heyting algebra $\widehat{\mathscr{C}}$ is isomorphic to $H$;
$R/H_{0}$ is atomless or has 2 elements depending on whether $\varepsilon$
takes the value $0$ or $1$; and $R/J$ is atomless for all ideals
$J\subsetneqq H_{0}$. He showed further~\cite[Corollary~2.2]{Alaev}
that any $\omega$-categorical $I$-algebra $R$ is isomorphic to
a product $R_{1}\times R_{2}\times\cdots\times R_{k}$, where $R_{k}$
has type $(H_{k},\varepsilon_{k})$ and $(H_{k},\varepsilon_{k})$
is as above, and that a minimal decomposition of this form is unique.
He also proved that an $\omega$-categorical $I$-algebra is uniquely
determined by $\widehat{\mathscr{C}}$ together with information about
the size and number of atoms of $R/J$ for each $J\in\widehat{\mathscr{C}}$
(see Remark~\ref{rem:Alaev}).

The pair $(H,\varepsilon)$ corresponds to a PO system $P$ with a
unique minimum element~$p$, such that $p\in P^{d}$ iff $\varepsilon=1$,
where $H$ is isomorphic to the Heyting algebra of upper subsets of
$P$. The minimum decomposition $R_{1}\times R_{2}\times\cdots\times R_{k}$
corresponds to a decomposition of $X$ into disjoint trim sets (see
Proposition~\ref{Propn std split of trim set}\ref{enu:Prop1.4}),
with $H_{k}=\{q\in P\mid q\geqslant p_{k}\}$ for some $p_{k}\in P_{\min}$
and $\varepsilon_{k}=1$ iff $p_{k}\in P^{d}$. We note the usual
order reversal between ideals (open sets) in $H(R)$ and closed sets
/ partition elements in the Stone space.

Underpinning the above discussion are the two parallel sets of dualities
mentioned in the introduction:
\begin{enumerate}
\item The dualities between finite closure algebras generated by their closed
elements, finite posets, and finite Heyting algebras;
\item The dualities between finite TBAs generated by their closed elements,
finite PO systems, and finite Heyting sa-algebras.
\end{enumerate}

\end{document}